\documentclass[12pt]{amsart}
\usepackage{amscd}
\usepackage{verbatim}
\usepackage{amssymb, amsmath, amsthm, amscd,ifthen}
\usepackage[dvips]{graphics}
\usepackage[cp866]{inputenc}
\usepackage{graphicx}
\usepackage{epsfig}
\usepackage{xcolor}

\usepackage{amsmath,amssymb,amscd,}
\usepackage[all]{xy}
\usepackage{url}

\usepackage{amssymb}

\newtheorem{thm}{Theorem}[section]

\newtheorem{prop}[thm]{Proposition}

\newtheorem{Ex}[thm]{Example}

\newtheorem{lemma}[thm]{Lemma}

\theoremstyle{definition}
\newtheorem{rem}[thm]{Remark}
\newtheorem{dfn}[thm]{Definition}

\usepackage{tikz}
\usepackage{tkz-euclide}
\usetikzlibrary{shapes,arrows,spy,positioning,snakes}

\title{Cooperative envy-free division}

\author[D. Joji\'{c}]{Du\v{s}ko Joji\'{c}}
\author[G. Panina]{Gaiane Panina}
\author[R. \v{Z}ivaljevi\'{c}]{Rade \v{Z}ivaljevi\'{c}}

\address[D. Joji\'{c}]{ Faculty of Science, University of Banja Luka}
\address[G. Panina]
{ St.\ Petersburg State University; St.\ Petersburg Department of Steklov Mathematical Institute}
\address[R. \v{Z}ivaljevi\'{c}]{Mathematical Institute of the Serbian Academy of Sciences and Arts (SASA), Belgrade}

\address{
  }

 \keywords{envy-free division, configuration space/test map scheme}

\begin{document}

\begin{abstract}

Relying  on configuration spaces and equivariant topology, we study a general \textit{``cooperative envy-free division problem''} where the players have more freedom of expressing their preferences (compared to the classical setting of the Stromquist-Woodall-Gale theorem).

A group of  players want to cut a ``cake'' $I=[0,1]$ and divide among themselves the pieces in an envy-free manner. Once the cake is cut and served in plates on a round  table (at most one piece per plate), each player makes her choice by pointing at one (or several) plates she prefers.
The novelty is that her choice may depend on the whole \textit{allocation configuration}. In particular, a player may choose an empty plate (possibly preferring one of the empty plates over the other), and take into account not only the content of her preferred plate, but also the content of the neighbouring plates.

We show that if the number of players is a prime power, in this setting an envy-free division still exists under standard assumptions that the preferences are closed.
\end{abstract}

 \maketitle \setcounter{section}{0}

\section{Overview and an informal introduction}\label{sec:Overview}

Extending and developing further the ideas from our earlier papers  \cite{jpz-2,PZ1,PZ2},  we  describe a new framework for applying  methods of equivariant topology to problems of mathematical economics, related to fair and envy-free division.

  In the classical approach,  exemplified by Stromquist-Woodall-Gale theorem, \cite{Strom, Wood, G}
the ``cake'' $I=[0,1]$  is cut into $r$ tiles, each player makes her choice by pointing at one (or several) tiles. The standard assumptions are that  (1) \textit{the preferences are closed}, and  (2) \textit{nobody prefers a degenerate tile}. Then
an envy-free division  exists for any number of players.

\medskip Offering a fresh perspective on the problem, we build a new narrative, where the players have more sophisticated preferences, may choose degenerate tiles (more precisely,  they may choose empty plates),   and the proofs rely on new \emph{configuration spaces and equivariant topology}.

\subsection{Players at a round table and a cake}
Assume that there are $r$ guests at a birthday party, who want to divide a birthday cake in such a manner that everyone is pleased with their piece of the cake.
(More generally, a group of $r$ players divide a commodity, also referred to as the ``cake'', modeled as the interval $I=[0, 1]$.)

\smallskip

Extending the narrative, we enrich this model by assuming that there is a round table with $r$ indistinguishable plates on it, so both the table and the set $\Pi$ of  plates have  the symmetry of the cyclic group $\mathbb{Z}_r$.

The cake is cut in at most $r$ pieces  (that is, in at most $r$ non-degenerate segments, referred as the ``tiles''), and the tiles are allocated to plates, at most one tile per plate. In some cases there might be empty plates.
We also do not exclude the case when there is only one tile  (the entire cake), and all the plates except one are empty.


Once a cut and an allocation is fixed, each player makes his choice by pointing at one of the plates he prefers (or more than one if they are equally valued and desired).
The novelty is that the choice may now depend on the \textit{allocation configuration} (allocation function) as a whole. Recall that in usual models the player would focus on a single (most desired) tile and completely ignore where the remaining (less desired) tiles are positioned.

\medskip
The following ``chocolate example'' illustrates a situation when the preferences depend on the cyclic ordering of the tiles on the table.

\bigskip
\noindent
\textbf{A ``chocolate example'' for a prime number of players, $r=p$.} Assume that the cake $I = [0,1]$ is non-uniformly decorated by marzipan, chocolate, cream, and other tasty ingredients.

Imagine a player wishes to maximize the total amount of chocolate
both on his chosen piece of cake and the two neighboring pieces, taken together.

 As a consequence it may happen (Example \ref{ex:empty-plate}), that the chosen plate is actually empty. Moreover, a player may prefer one empty plate over another, i.e\ the empty plates are not necessarily equal!

Observe that the preferences do not change if we rearrange the pieces of the cake by cyclically permuting  the plates on the table.

\medskip

The rationale behind the choice of an empty plate may be that a player anticipates a future \textit{cooperation} or \textit{trade} with his immediate  neighbors at the table.  Some other players may have similar plan, regarding other ingredient of the cake, say the cream or the marzipan.

\medskip
A special case of our Theorem \ref{thm:main-gen}  states that in this case the envy-free division is always possible, under  some quite natural and mild assumptions on the preferences.


\medskip\noindent
\textbf{A ``chocolate example'' for a prime squared, $r=p^2$.} Among the corollaries of Theorem \ref{thm:main-gen} is the following result, addressing the case when the number of players is a square of a prime.

\medskip
As in the previous example we  assume that the cake $I$ is (non-uniformly) decorated by marzipan, chocolate, cream, etc. This time
there are $p$  round tables, each having $p$ plates (so the cardinality of the set $\Pi$ of all plates is $r=p^2$). The tiles of the cake are allocated to the plates, and a player again may wish to maximize the total amount of chocolate
on his own piece and the two neighboring pieces taken together.  The conclusion is that, no matter what the other players prefer, the envy-free division always exists.

\medskip
In the general case the players may have reasonably complex, unknown or even ``irrational'' objectives and preferences, from an outside point of view.
What really matters is that the preferences do not change if we rearrange the pieces of the cake by cyclically permuting the tables and (again cyclically) permute the plates at each table.

\medskip
\noindent
\textbf{A variation on the theme of ``chocolate examples''.}  Imagine   that for each player $j\in [r]$ there is a  continuous \textit{evaluation function} $f^j$, defined on the collection  of all segments $[a,b]\subseteq [0,1]$, specifying the value of each segment from the viewpoint of that player. The only condition is that $f^j$  should be constant on all degenerate segments (that is in the case $a=b$).

Assume that the plates a positioned on the table in a circular order. The player  $j$  chooses a plate number $i$ provided  $$3f^j(\hbox{segment  lying in the plate }i)-f^j(\hbox{segment  lying in the plate }i+1) + $$$$f^j(\hbox{segment  lying in the plate }i-1))$$ attains its maximum.
(If the maximum is attained at two more or positions, he chooses each of them.) Then, as a consequence of Theorem \ref{thm:main-gen}, an envy-free division always exists if $r=p$ is a prime number.

\subsection{Novelty of the model}

The simple ``chocolate examples'' illustrate some of the main features of the proposed new approach to the ``cooperative envy-free division''. The novelty of the methods and ideas used in the paper might be, for the convenience of a casual or non-expert reader, summarized as follows.

\begin{enumerate}
  \item The preference of a player depends on the whole allocation configuration. In other words the player takes into account the entire allocation function $\alpha : [k] \rightarrow \Pi$ (where $k\leq r$), which describes the placement of the pieces in the plates. Recall that in standard models a player would focus on a single (most desired) tile and completely ignore where the remaining (less desired) tiles are positioned.
  \item If empty plates are offered, which happens if the number of pieces of the cake is strictly smaller than $r$, then the players are allowed to choose one of the empty plates. Moreover, an important novelty is that the empty plates are not necessarily equal, i.e.\ one empty plate may be preferred over the other.
\end{enumerate}

\begin{Ex}\label{ex:empty-plate}
One can easily build a ``chocolate example'', with two empty plates, illustrating {\rm (2)}. This happens, for instance, with seven plates having the following distribution of the chocolate $$(100,0,100,1,1,0,1)$$ where a  player prefers   one of the empty plates, but not the other.
\end{Ex}

For comparison,  in the classical Stromquist-Woodall-Gale theorem, \cite{Strom, Wood, G}, the players are not allowed, under any circumstances, to choose empty plates.

In the more recent approach of Avvakumov and Karasev  \cite{AK}, the players are allowed to choose an empty plate. However, in their approach all empty plates are considered equal and equally desired (or undesired).

\begin{enumerate}
  \item[\rm (3)] A technical novelty is the use of a  new configuration space $\mathcal{C}_r$ (Section \ref{sec:form}), as a natural mathematical model for a cooperative envy-free division. Recall that in usual applications of topological methods, see for example \cite{Strom, Wood, G, AK, Fr, MeZe, S-H},  the model for the space of cuts of the cake into $r$ pieces is the $(r-1)$-dimensional simplex $\Delta^{r-1}$.
\end{enumerate}

\subsection{Summary of the new results}

What is the natural environment (ecological niche) for all our ``chocolate examples'', and other similar ``cooperative envy-free division results''?
What is the role of symmetries of preferences, which ``do not change if we cyclically rearrange the pieces of the cake'', etc.? Is the role of prime numbers essential (intrinsic) or they appear as an artefact, reflecting the limitations of the methods applied.

\medskip
Our main new result (Theorem \ref{thm:main-gen}, Section \ref{sec:statements}) addresses all these questions.
It guarantees the existence of a cooperative envy-free division under the following conditions: \textit{the number $r$ of players is a prime power},  \textit{the preferences are closed}, and the \textit{preferences  respect the action of the $p$-toral group on the set $\Pi$ of plates}.

\medskip

  Both the formulation and the proof of Theorem \ref{thm:main-gen} reveal that the symmetries of preferences are indeed ``at the heart of the matter''.

In particular it is known that in closely related problems,  see \cite[Section 4.3]{AK} and \cite[Theorem 4.4]{PZ1},  the order of the group of symmetries must be a prime power. By a similar argument this condition is also unavoidable in the ``cooperative envy-free division problem''.

\medskip
The proof of Theorem \ref{thm:main-gen} is postponed for Section \ref{sec:proofs}, after the introduction of relevant topological tools and ideas, see also the Appendix (Section \ref{sec:appendix}) for some auxiliary material and a guide to the literature.

\medskip

Theorem \ref{thm:main-gen} paves the way for other envy-free division results for cooperative preferences. This holds, in particular, for the envy-free division in the presence of a ``dragon'' \cite{PZ2}, where we distinguish two basic scenarios:

\medskip

(1)  There are $r-1$ players and a dragon. Once the ``cake'' is divided into $r$ pieces the dragon, ignoring all other players, makes his choice and grabs one of the pieces. After that the players want to divide the remaining pieces in an envy-free fashion.

(2) There are $r+1$ players who divide the cake into $r$ pieces. A ferocious dragon comes and swallows one of the players. The players want to cut the cake in advance in such a way that no matter who is the unlucky player swallowed by the dragon, the remaining players can share the tiles in an envy-free manner.
\bigskip

We show (Theorems \ref{thm:(r-1)} and \ref{thm:(r+1)}  in Section \ref{sec:statements}) that in both of these scenarios one can guarantee the existence of envy-free division for cooperative preferences, under similar conditions as in Theorem \ref{thm:main-gen}.

\section{Mathematical formalization of the problem via configuration spaces and cooperative preferences}\label{sec:form}

For $r\in \mathbb{N}$, let us describe a topological space  $\mathcal{C}=\mathcal{C}_r$ whose points correspond to cuts of the cake and allocations of the tiles to $r$ plates.
The space will be called the \textit{natural configuration space}, or simply the ``configuration space'', for the problem of cooperative envy-free division.

We call it ``natural'' to emphasize that it records, in simplest mathematical terms, all possible scenarios of cutting the cake and distributing (allocating) the pieces to the plates (at most one piece per plate).

\medskip
\begin{itemize}
  \item
\textit A {proper  cut} $x$ (or cut, for short) of length (tile number) $k$ is a sequence
$$0=x_0<x_1<x_2<...<x_{k-1}<x_k =1, $$ where $k\leq r$.  Each proper cut creates a partition  of $I$ into $k$ non-degenerate segments $I_i(x) = [x_{i-1}, x_i]$,  called the \textit{tiles}, which are numbered in the increasing order, from left to the right.
(We do not exclude the possibility of an ``empty cut'' with the only one tile, corresponding to the case when $k=1$.)
  \item
For a fixed cut $x$ with the tile number $k$, an \textit{allocation function}  is an injective function $$\alpha: [k]\rightarrow \Pi$$
where $\Pi$ is the set of ``plates'', arranged at a round table. For example one can imagine that there are $r$ plates numbered by $1,...,r$ and that the tile $I_i = I_i(x)$  is served in the plate $\alpha(i)$.  Note that some of the plates may remain empty.
Injectivity means that the tiles are allocated to the plates, with the condition \textit{``at most one  tile per plate''}.

\medskip
We define the \textit{configuration space} $\mathcal{C}=\mathcal{C}_r$ as the set whose elements are all possible pairs
$$(x,\alpha) \,  =  \, \mbox{(a proper cut, an allocation function)} \, . $$
\end{itemize}

\subsection{Topology on the configuration space $\mathcal{C}$}
\label{sec:topo}

Topology on $\mathcal{C}$  is defined in two natural and equivalent ways.

\medskip\noindent
(A) (Topology via converging sequences)

\medskip
  A sequence $(x^n,\alpha_n)_{n=1}^{\infty}$ converges to $(x,\alpha)$ if and only if: \medskip

\begin{enumerate}
  \item The sequence  of  proper cuts $x^n =  (x_0^n,\dots ,x_{k_n}^n)\!$  converges to the proper cut $x = (x_0,\dots ,x_{k})$ in the sense of Hausdorff metric, see the
  Appendix (Section \ref{sec:symm-power}). In light of the ``stabilization property'' (Proposition \ref{prop:stabilization}) this means that
  for each $n$ there is a subcollection $\mathcal{T}_{n} = \{J^n_i\}_{i=1}^k \subseteq \mathcal{T}_{x^n}$ (of \emph{``essential tiles''} of $x^n$) such that the associated sequences of intervals converge $J^n_i \longrightarrow I_i$ to the corresponding tiles in $\mathcal{T}_x$, while the lengths of all other tiles from $\mathcal{T}_{x^n}\setminus \mathcal{T}_{n}$ tend to zero when $n \rightarrow +\infty$.

  \item The values of $\alpha_n $ ``stabilize'' on essential tiles $\mathcal{T}_{n} = \{J^n_i\}_{i=1}^k$ and coincide with the corresponding values of the allocation function $\alpha$. More precisely, there exists $n_0$ such that for each $n\geqslant n_0$ and each $i=1,\dots, k$
      \[
          \alpha_n(i) =: \alpha_n(J_i^n) = \alpha(I_i) := \alpha(i) \, .
      \]
\end{enumerate}
Note that in (1) the sequence $(k_n)_{n=1}^{+\infty}$ of ``tile numbers'' is not necessarily convergent (let alone convergent to $k$). Note also that in (2) the values of the allocation function $\alpha_n$ on inessential tiles $\mathcal{T}_{x^n}\setminus \mathcal{T}_{n}$ is not relevant at all for the convergence of $(x^n,\alpha_n)_{n=1}^{\infty}$ to $(x,\alpha)$.

\medskip\noindent
(B) (Topology  via neighborhoods)

\medskip
Given a proper cut  and an allocation function $(x,\alpha)$, let us describe a neighborhood $O_\varepsilon (x,\alpha)$  of this point in the configuration space $\mathcal{C}_r$. Set $\varepsilon$ to be much smaller than the length of the  minimal tile.

Then $(y,\beta)\in O_\varepsilon (x,\alpha)$ if and only if
\begin{enumerate}
  \item  The Hausdorff  distance between $x$ and $y$ is less than $\varepsilon$;
  \item  If $J\in \mathcal{T}_{y}$ such that there exists (necessarily unique) tile $I\in \mathcal{T}_{x}$ such that the length of $I\cap J$ is much bigger than $\epsilon$, then $\beta(J) = \alpha(I)$.
\end{enumerate}

\subsection{Group action on $\mathcal{C}_r$}
\label{sec:G-action}

From here on we assume that a group $G$ acts (from the left) on the set $\Pi$ of plates, as a group of permutations.  By default this action is free and transitive. This induces a free action on the configuration space $\mathcal{C}_r$ where for $(x,\alpha)\in \mathcal{C}_r$,
\[
     \sigma (x,\alpha):= (x, \sigma \circ\alpha) \, .
\]
At this stage we usually assume that $\Pi = [r] = \{1,\dots, r\}$ where $r=p^\nu$ is a prime power and $G = (\mathbb{Z}_p)^\nu$ is a $p$-toral group. In particular, for the first chocolate example, $G = \mathbb{Z}_p$  is the cyclic group permuting cyclically the plates on the round table.

\subsection{Preferences}\label{sec:A-preferences}


The natural configuration space $\mathcal{C}_r$ is the ambient space for \emph{cooperative preferences} $A^j_i$ where $j\in [r]$ is a label of a player and $i\in \Pi$ is a label of a plate. (In concrete applications usually $\Pi = [r] = [p^\nu]$.)

\begin{dfn}\label{dfn:A-pref}
The preferences of $r$ players is a collection (matrix) of subsets $(A^j_i)_{i\in \Pi}^{j\in [r]}$ of the configuration space $\mathcal{C}_r$, indexed  by $(i,j)\in \Pi \times [r]$. The subsets are interpreted as preferences as follows:
\begin{equation}\label{eqn:prefs-2}
\begin{split}
f=(x, \alpha) \in A^j_i \ \ \Leftrightarrow  \, &  \mbox{ {\rm in the cut} } x \mbox{ {\rm and the allocation} }
\alpha \mbox{ {\rm the player} } j\\
& \mbox{ {\rm prefers $\alpha^{-1}(i)$, the content of the plate } } i\in \Pi \, .
\end{split}
\end{equation}
\end{dfn}

\begin{dfn}\label{def:default} By default all preferences are \emph{closed, covering} and \emph{equivariant}.
\begin{enumerate}
  \item[\rm (a)]  The preferences are closed if $A^j_i$ are closed subsets of $\mathcal{C}_r$.
  \item[\rm (b)] The preferences are covering if\, $\bigcup_{i\in \Pi} A^j_i = \mathcal{C}_r$ for each $j\in [r]$.
  \item[\rm (c)] The preferences are equivariant if  for each $\sigma\in G$,
\begin{equation}\label{eqn:equivariance-1}
  (x, \alpha ) \in A_{i}^j  \Leftrightarrow  \sigma(x, \alpha )\in A_{\sigma(i)}^j   \, .
\end{equation}
\end{enumerate}
\end{dfn}
Let us take a closer look at the condition (\ref{eqn:equivariance-1}). As a consequence of (\ref{eqn:prefs-2}),  $(x, \sigma\circ \alpha) \in A^j_{\sigma(i)}$ if and only if in the cut $x$ and the allocation $\sigma\circ \alpha$ the player $j$
prefers the content of the box $\sigma(i)$. Since $(\sigma\circ\alpha)^{-1}(\sigma(i)) = \alpha^{-1}(i)$ the condition (\ref{eqn:equivariance-1}) expresses the idea of ``invariance of preferences'' with respect to the group action. Informally, the players
make their decisions on the basis of the ``inner structure'' of the whole allocation function and the equivariance condition guarantees that this inner structure is preserved by the group action.

\begin{Ex}{\rm
 In  our ``chocolate examples'', reviewed in the Introduction, the set $\Pi$ was interpreted as a set of plates distributed around circular tables. The invariance of preferences with respect to the circular action of the group $G$ on $\Pi$ has a clear geometric meaning, clarifying what is in this case meant by the ``inner structure'' of an allocation function.  A more general interpretation would involve an abstract ``table'' $\Pi$ with an associated $G$-invariant graph (hypergraph) on $\Pi$, encoding possible patterns of cooperation between players.  }
\end{Ex}

The following definition formalizes, in the language of the matrix of preferences, the precise meaning of the cooperative envy-free division.

\begin{dfn}\label{def:coop} Let $(A^j_i)_{i\in \Pi}^{j\in [r]}$ be a matrix of  preferences.
A point $  (x,\alpha) \in \mathcal{C}$  provides a \emph{ cooperative envy-free division} for preferences  $(A^j_i)_{i\in \Pi}^{j\in [r]}$ if  there exists a bijection $\sigma : [r] \rightarrow \Pi$ such that
\begin{equation}\label{eqn:non-emty-envy-free}
(x,\alpha)\in \bigcap_{j=1}^r A_{\sigma(j)}^j  \, .
\end{equation}
\end{dfn}

\section{Statements of the main results }\label{sec:statements}

Our main result is the following ``cooperative envy-free division theorem''.

\begin{thm}\label{thm:main-gen}
  Suppose that $r=p^\nu$ is a prime power. Let $G = (\mathbb{Z}_p)^\nu$ be a $p$-toral group which acts on the set of ``plates'' $\Pi \cong [r]$ by a free and transitive action. Let $(A^j_i)^{j\in [r]}_{i\in \Pi}$ be a matrix of subsets of the \emph{natural configuration space} $\mathcal{C}_r$, defining the preferences of the players. Assume that these preferences  are closed, covering and equivariant in the sense of Definition  \ref{def:default}. Then there exists a  cooperative envy-free division of the cake satisfying the preferences of all players (Definition \ref{def:coop}).
\end{thm}

\medskip

The following two theorems extend Theorem \ref{thm:main-gen} to the case when one of the players, called the dragon, is greedy and non-cooperative, which makes the envy-free division even more difficult.

\medskip
There are two basic scenarios. In the first scenario after the cake is divided into at most $r$ pieces and distributed in $r$  plates $\Pi$ (at most one piece per plate) the dragon, ignoring all other players,  grabs one of the plates. After that the remaining players should be given the rest of the plates so that there is no envy between them  (dragon included!).

\begin{thm}(The dragon takes a piece of the cake) \label{thm:(r-1)}
Let $r$ be a prime power. Let
  $$
  (A_i^j)_{i\in\Pi}^{j\in [r-1]},\ \  A_i^j\subseteq \mathcal{C}_r
  $$
be a $\Pi\times (r-1)$-matrix of preferences which are closed, covering, and equivariant. Then one can choose for each $j\in [r-1]$ two distinct elements $u_j$ and $v_j$ in $\Pi$ such that
 \begin{enumerate}
   \item[(a)] The collection $E =\{e_j\}_{j\in [r-1]}$ of two element sets $e_j=\{u_j, v_j\}$ is the edge-set of a tree $T=(V,E)$  on $V = [r]$.
   \item[(b)]    \begin{equation}\label{eqn:tree-intersection-2}
                 \bigcap_{v \mbox{ {\rm\tiny incident to} }e } A_v^e  =  \bigcap_{j\in [r-1]}   (A_{u_j}^j\cap A_{v_j}^j) \neq \emptyset \, .
   \end{equation}
 \end{enumerate}
\end{thm}

In the second scenario there are $r+1$ players who divide the cake into $r$ pieces. A ferocious dragon comes and swallows one of the players. The players want to cut the cake in advance, and put them in plates in such a way that no matter who is the unlucky player swallowed by the dragon, the remaining players can share the pieces in an envy-free manner.
The following theorem says that is always possible, under the same assumptions on the preferences as in Theorem \ref{thm:main-gen}.

\begin{thm}(The dragon takes a player)  \label{thm:(r+1)}
Let $r$ be a prime power.
  Let
  $$
  (A_i^j)_{i\in\Pi}^{j\in [r+1]},\ \  A_i^j\subseteq \mathcal{C}_r
  $$
be a $\Pi\times (r+1)$-matrix of preferences which are closed, covering, and equivariant.

 Then one can choose for each $i\in\Pi$ two distinct elements $u_i$ and $v_i$ in $[r+1]$ such that
 \begin{enumerate}
   \item The collection $E =\{e_i\}_{i\in \pi}$ of two element sets $e_i=\{u_i, v_i\}$ is the edge-set of a tree $T=(V,E)$  on vertices $V = [r+1]$.
   \item    \begin{equation}\label{eqn:tree-intersection-333}
                 \bigcap_{v \mbox{ {\rm\tiny incident to} }e } A_e^v  =  \bigcap_{i\in [r]}   (A^{u_i}_i\cap A^{v_i}_i) \neq \emptyset \, .
   \end{equation}
 \end{enumerate}
\end{thm}

\section{Proofs via new configuration spaces, equivariant topology, chessboard complexes and the Birkhoff polytope}
\label{sec:proofs}

Following the main idea of \cite{PZ1}  (which originally evolved from \cite{jpz-2}),
we  change the basic setup and modify the narrative, by allowing more cut points, which are allowed to coincide and create degenerate tiles.

\medskip
More explicitly, an \emph{improper cut} of the segment $I=[0,1]$ by $n-1$ cut-points is  recorded  as  $x=(0\leqslant x_1\leqslant \dots \leqslant x_{n-1} \leqslant 1)$. The tiles $I_i = I_i(x) = [x_{i-1}, x_i]$ are automatically labeled by $i\in [n]$,  from left to the right. The novelty, compared to Section \ref{sec:form}, is that the inequalities in the cut are no longer strict and, as a consequence, some of the tiles $I_i$ may be degenerate.

\medskip
The introduction of degenerate tiles (albeit unnatural from the view point of cutting a real cake) has the advantage that the
\emph{configuration space of all improper cuts} is the standard $(n-1)$-dimensional simplex $\Delta^{n-1}$.

\medskip
Recall that in the more realistic model, where we used the space $I(r)$  of all proper cuts of tile number at most $r$ (Sections \ref{sec:form} and \ref{sec:symm-power}), the degenerate tiles are not permitted. However, it is well known \cite{ams-93, ST} that the topological structure  of symmetric powers $I(r)$ is in general much more complex and less transparent than the simple geometry of a simplex.

\medskip
This observation motivates us to replace \emph{proper} by \emph{improper} cuts in the definition of the natural configuration space $\mathcal{C}_r$, introduced in
Section \ref{sec:form}, hoping that the new ``auxiliary configuration space'' $\mathcal{C}'_r$ is simpler and easier to handle.

\subsection{Auxiliary configuration space $\mathcal{C}'_r$}

\begin{dfn}  Given an improper cut $x=(0\leqslant x_1\leqslant \dots \leqslant x_{n-1} \leqslant 1)$ of the segment into $n$ tiles $\{I_i(x)\}_{i=1}^n$, an \textit{allocation function}  of the tiles into a set $\Pi$ of plates  is an ``essentially injective'' map
$$\alpha:[n]\rightarrow \Pi \, .$$
Essential injectivity of $\alpha$ means that \textit{not more than one non-degenerate tile
is  allocated to a single plate}. In particular, if $n>r := \vert \Pi\vert$, not every cut allows an allocation function.

\medskip
A pair $(x,\alpha)$ is called an (improper) partition/allocation of the ``cake'' $[0,1]$.
Two such pairs $(x,\alpha)$ and $(x,\alpha')$, with the same improper cut $x$, are called \emph{equivalent} if they differ only by allocations of degenerate tiles. We write $[(x,\alpha)]$ for the corresponding equivalence class.
\end{dfn}

\begin{rem}
  The last part of the definition formalizes (on the level of configuration spaces) the property
``all degenerate tiles are equal'' having the ``zero value''.
\end{rem}

\begin{Ex}
Let $n=6$ and assume  that the cut-points (with repetitions) are $(1/4,1/3,1/3,2/3,2/3)$. There are altogether six tiles $\{I_i\}_{i=1}^6$ and two of them, the tiles $I_3$ and $I_5$, are degenerate. Let $\Pi = \{\Pi_i\}_{i=1}^6$ be a set of six plates.

Then the following three allocations are equivalent:
\newline $\Pi_1=\{1\}, \Pi_2=\{6\}, \Pi_3=\{5\}, \Pi_4=\{4\}, \Pi_5=\{3\}, \Pi_6=\{2\}$
\newline $\Pi_1=\{1\}, \Pi_2=\{6\}, \Pi_3=\{3\}, \Pi_4=\{4\}, \Pi_5=\{5\}, \Pi_6=\{2\}$
\newline $\Pi_1=\{1,3,5\}, \Pi_2=\{6\}, \Pi_3=\emptyset, \Pi_4=\{4\}, \Pi_5=\emptyset, \Pi_6=\{2\}$.

\medskip
It is instructive to think that for an equivalence class $[(x,\alpha)]$ only non-degenerate tiles are allocated, whereas degenerate ones are ignored (and put in the ``trash'').
So the equivalence class of the above  partition/allocation is recorded as
\newline  $\Pi_1=\{1\}, \Pi_2=\{6\}, \Pi_3=\emptyset, \Pi_4=\{4\}, \Pi_5=\emptyset, \Pi_6=\{2\}$.
\end{Ex}

\begin{dfn}\label{dfn:auxiliary-conf-space}
  Assume that the set of plates is $\Pi = \{\Pi_i\}_{i=1}^r \cong [r]$ and choose $n=2r-1$. Define the auxiliary configuration space $\mathcal{C'} = \mathcal{C}'_r = \{[(x,\alpha)]\}$ as the set of all (equivalence classes of) partitions/alocations $(x,\alpha)$, where $x\in \Delta^{n-1}$ is an improper partition (cut) of $[0,1]$ and $\alpha : [n] \rightarrow \Pi$ a corresponding allocation function.
\end{dfn}
$\mathcal{C}'_r$ is, according to Definition \ref{dfn:auxiliary-conf-space}, just a set. However, it is not difficult to see that it has a natural topology. Indeed, for a fixed allocation function $\alpha : [n] \rightarrow [r]$ the set $\Delta_\alpha \subset \mathcal{C}'_r$ of all classes $[(x,\beta)]\in \mathcal{C}'_r$ such that $\beta = \alpha$
is naturally isomorphic to the simplex $\Delta^{n-1}$.

\medskip
Therefore $\mathcal{C}'_r = \cup_\alpha \Delta_\alpha$ is just the union of these simplices, glued together. In order to see this gluing more clearly, we associate to the partition/allocation $(x,\alpha)$ a $(r\times n)$-matrix $M_{(x,\alpha)}$ with real entries
\begin{equation}\label{eqn:matrices}
   M_{(x,\alpha)} =  (x_1 - x_0)E_{1,\alpha(1)} +(x_2 - x_1)E_{2,\alpha(2)} + \dots + (x_{n} - x_{n-1})E_{n,\alpha(n)}
\end{equation}
where $E_{i,j}$ is the $(r\times n)$-matrix having zeros everywhere, except at the position $(i,j)$ where the corresponding entry is $1$.

\medskip
It is not difficult to see that $M_{(x,\alpha)}$ depends only on the equivalence class of $(x,\alpha)$ and that $M_{(x,\alpha)} = M_{(y,\beta)}$ if and only if $y=x$ and $(x,\alpha)$ and $(x,\beta)$ are equivalent.

\medskip
It follows from the faithfulness of representation (\ref{eqn:matrices}) that $\mathcal{C}_r'$ is homeomorphic to (the geometric realization of) a simplicial complex, piecewise linearly embedded in the linear space $Mat_{r\times n}^\mathbb{R}$ of all $(r\times n)$-matrices. { It turns out that this simplicial complex is a well-studied object in topological combinatorics, known under the name ``chessboard complex'' \cite{blvz, Z17}, see also Section \ref{sec:chess} for a short overview.}

The following proposition is of fundamental importance in the theory of chessboard complexes.

\begin{prop}\label{prop:chessboard}
The configuration space  $\mathcal{C}'_r$ is isomorphic to the geometric realization of the chessboard complex $\Delta_{r,2r-1}$.
\end{prop}
\proof
Indeed, the chessboard complex $\Delta_{r,n}$ is described \cite{Z17} as the complex of all non-attacking configurations of rooks on a $(r\times n)$-chessboard. The graph $\Gamma(\alpha)\subset [n]\times [r]$ of each allocation function $\alpha : [n]\rightarrow [r] $ describes such a configuration of rooks. Conversely, each maximal simplex (arrangement of rooks) in $\Delta_{r,n}$ arises uniquely from such an allocation function. \qed

\medskip
A permutation group $G$, acting on the set $\Pi \cong [r]$ of plates,  acts on the natural configuration space $\mathcal{C}_r$. A similar action exists on the auxiliary configuration space $\mathcal{C}_r'$ where $\sigma(x,\alpha)= (x, \sigma \circ \alpha)$ for each $\sigma\in G$. From the view point of the complex $\Delta_{r,2r-1}$, this is the
action arising from permuting the rows of the chessboard $[r]\times [n]$.


\bigskip
{ The auxiliary configuration space $\mathcal{C}'_r \cong \Delta_{r,2r-1}$ is introduced to mimic, for technical purposes, the natural configuration space $\mathcal{C}_r$. For the same reason we introduce {special sets} (``ghost preferences'') $B^j_i\subseteq \mathcal{C}'_r$, mimicking the role of actual preferences $A^j_i\subseteq \mathcal{C}_r$, as arbitrary subsets of $\mathcal{C}'_r$ which are closed, covering and equivariant in the sense of Definition \ref{def:default}. }

\medskip
The ultimate justification for introducing the auxiliary configuration space $\mathcal{C}'_r$ is the following result, originally proved in \cite{PZ1}.

\begin{thm}\label{ThmAux}
Let $r = p^\nu$ be a prime power.
  Let $$(B_i^j)_{i,j=1}^r,\ \  B_i^j\subseteq \mathcal{C}'_r$$ be a matrix of  closed, covering sets which are equivariant with respect to the group $G=(\mathbb{Z}_p)^\nu$. Then there exists a bijection $\sigma : [r]\rightarrow [r] \cong \Pi$ such that
\begin{equation}\label{eqn:empty-4}
     \bigcap_{j\in [r]}  B_{\sigma(j)}^j \neq\emptyset \, .
\end{equation}
\end{thm}
\proof

Although a detailed proof can be found in \cite{PZ1} we repeat, for the sake of completeness, some of the central steps. Moreover, an alert reader will see (at the end of the proof) where exactly the use of the auxiliary configuration space $\mathcal{C}'_r \cong \Delta_{r,2r-1}$ was important.

Assume, for the sake of contradiction, that the intersection (\ref{eqn:empty-4}) is empty.
Construct an (equivariant) test map $F : \mathcal{C}' \rightarrow \mathbb{R}^r$ which records the information provided by the preferences $B_i^j$.

Initially we replace {the collection of sets} $B_i^j $ by {a collection of slightly larger open sets}  $O_i^j$ and choose  (for each $j\in [r]$) an equivariant partition of unity $\{f_i^j\}_{i=1}^r$ subordinated to the cover $\{O_i^j\}_{i=1}^r$.
\medskip
Since the preferences are equivariant, we can  make these  functions also equivariant:
\begin{equation}\label{eqn:equivariance-3}
f_{\sigma(i)}^j(\sigma(x, \alpha)) = f_{i}^j(x, \alpha) \, .
\end{equation}
By averaging $F_i=\frac{1}{r}\sum_{j=1}^{r}f^j_i$ we obtain a vector-valued $S_r$-equivariant function  $$F=(F_1,F_2,\dots,F_r):\mathcal{C}' \rightarrow \mathbb{R}^r.$$
 Let $\widehat{F}  : \mathcal{C}' \rightarrow \mathbb{R}^r/D$ be the (equivariant) map obtained by composing the map $F$ with the projection $\mathbb{R}^r \rightarrow \mathbb{R}^r/D$, where $D = \{(y,\dots, y)\}_{y\in \mathbb{R}}\subset \mathbb{R}^r$ is the diagonal subset of $\mathbb{R}^r$. This map must have a zero.

 Indeed, in the opposite case there arises an equivariant map $\mathcal{C}' \rightarrow S^{r-2}$, which contradicts Volovikov's theorem \cite{Vol96-1}, in light of the fact that the chessboard complex $\Delta_{r,2r-1}$ is $(r-2)$-connected.

\medskip
The proof is completed by an application of the original idea of D. Gale, exploiting the fact the vertices of the Birkhoff polytope of bistochastic matrices are permutation matrices.
 \qed

\subsection{Proof of Theorem \ref{thm:main-gen}}\label{sec:pf-thm-ge-n}

There is a  natural  map  $$\pi:\mathcal{C}'\rightarrow \mathcal{C}$$
defined as follows. A  cut together an allocation of its non-degenerate tiles induces a proper cut (by forgetting degenerate tiles) with the  same allocation function.

\begin{lemma}
  The map $\pi$ is well-defined. Moreover, it is a non-injective epimorphism, which is both continuous and equivariant.
\end{lemma}
The proof is quite straightforward, for example the continuity follows by checking that for each sequence $z_n \rightarrow z$,  convergent in $\mathcal{C}'$, the sequence $\pi(z_n)$ converges to $\pi(z)$ in the sense of Section \ref{sec:topo}.

\medskip
Assuming the lemma the proof of Theorem \ref{thm:main-gen} is finished as follows.
Given the matrix of preferences $(A_i^j)$, where $A_i^j\subseteq \mathcal{C}$, let $B_i^j= (\pi)^{-1}(A_i^j)$ be the corresponding family of subsets of the auxiliary configuration space $\mathcal{C}'$.

By Theorem \ref{ThmAux}, there is an envy-free division for $(B_i^j)$, that is a bijection $\sigma : [r]\rightarrow  \Pi$ and an element
\[
(x,\alpha)\in \bigcap_{j\in [r]}  B_{\sigma(j)}^j \, .
\]
The image $\pi((x,\alpha))\in \mathcal{C}$ is clearly an envy-free division for $(A_i^j)$. \qed

\bigskip

The  dragon versions are routine in view of \cite{PZ2}.

\section{Appendix: Spaces of finite sets}\label{sec:appendix}


\subsection{Hausdorff metric} The Hausdorff metric $d_H(A,B)$ measures the distance between two (non-empty) closed sets $A$ and $B$ in a metric space $(X,d)$.

If $A = \{a\}$ is a point than (by definition)
\[
d_H(a,A):= d_H(\{a\},A):= \inf_{x\in A}\{d(a,x)\}
\]
and by  symmetry, $d_H(A,a) = d_H(a,A)$. In general,
\[
 d_H(A,B)= \max\{\sup_{a\in A} d(a,B), \sup_{b\in B} d(A,b) \}  \, .
\]
If $A$ and $B$ are finite (more generally compact) subsets of $X$ then, in these formulas, we are allowed to write  $\max$ instead of $\sup$ and $\min$ instead of $\inf$.

\medskip

If $(X,d)$ is the interval $[0,1]$ with the usual distance function $d(a,b) = \vert a-b\vert$, then $d_H(a,B)$ (respectively $d_H(A,b)$) is  the distance from $a$ to its nearest (left or right) neighbour in $B$ (similarly the nearest neighbour of $b$ in $A$), and $d_H(A,B)$ is simply the largest of all these distances.

\medskip
An alternative definition describes the Hausdorff distance $d_H(A,B)$ (for compact A and B) as the smallest $\epsilon \geqslant 0$ such that both $B\subseteq O_\epsilon(A)$ and
$A\subseteq O_\epsilon(B)$, where $O_\epsilon (A) = \{x\in X \mid (\exists a\in A)\, d(x,a) \leqslant \epsilon\}$ is the closed $\epsilon$-neighborhood of $A$ in $X$.

\subsection{Symmetric power}
\label{sec:symm-power}

 The $n^{th}$ \emph{symmetric power} of a metric space $X$ \cite{ams-93, ST}, denoted by $[X]^{\leqslant n}$, is the space of all non-empty subsets of $X$ of cardinality at most $n$. We assume that $[X]^{\leqslant n}$ is metrized with the Hausdorff metric. The following closely related space
$$I(r) :=  \{A \subset [0,1] \mid \{0,1\} \subseteq A \mbox{ {and} } \vert A\vert\leq r+1\} \, ,$$
where $I = [0,1]$ and $r\geq 1$, is also referred to as a symmetric power (of order $r$).

\medskip
Elements of $I(r)$ can be also described as strictly increasing sequences
\begin{equation}\label{eqn:seq}
 x : \quad 0 = x_0 < x_1 < \dots < x_{k-1} < x_k = 1
\end{equation}
of the length (tile number) $k$, where $1\leq k\leq r$. Such a sequence is often referred to as a \emph{proper cut} of the tile number $k$, where
$\mathcal{T}=\mathcal{T}_x = \{I_j(x)\}_{j=1}^k$ is the associated set of $x$-tiles, $I_j(x):= (x_{j-1}, x_j)$.

\medskip
The following elementary proposition records (for future reference) an important property of the Hausdorff convergence of cuts, which is referred to as the
``stabilization property'' of the corresponding sets of tiles.

\begin{prop}\label{prop:stabilization}
Let $x^n = (x_i^n)_{i=0}^{k_n}$ be a sequence of proper cuts $(n=1,2,\dots)$, with the associated sequences $(k_n)$ of tile numbers and sets $(\mathcal{T}_{x^n})$ of tiles. Suppose that
$x^n \longrightarrow x$ converges in Hausdorff topology to a cut $x = (x_i)_{i=0}^k$, with the tile number $k$ and the set of tiles $\mathcal{T}_x = \{I_i(x)\}_{i=1}^k$. Then $\mathcal{T}_{x^n}$ \emph{converges essentially} to $\mathcal{T}_{x}$
\[
        \mathcal{T}_{x^n} \stackrel{ess}{\longrightarrow} \mathcal{T}_{x}
\]
in the following sense. For each $n$ there is a subcollection $\mathcal{T}_{n} = \{J^n_i\}_{i=1}^k \subseteq \mathcal{T}_{x^n}$ (of \emph{``essential tiles''}) such that the associated sequences of intervals converge $J^n_i \longrightarrow I_i(x)$ to the corresponding tiles in $\mathcal{T}_x$. As a consequence the lengths of all other tiles from $\mathcal{T}_{x^n}\setminus \mathcal{T}_{n}$ tend to zero when $n \rightarrow +\infty$.
\end{prop}

\subsection{Topological configuration space $\mathcal{C}_r$}

\subsection{Chessboard complexes}\label{sec:chess}

\medskip
The \emph{chessboard complex} $\Delta_{m,n}$ is an abstract simplicial complex defined on an $m\times n$ chessboard with $m$ rows  and $n$ columns. More precisely,
the vertices of $\Delta_{m,n}$ are the squares of the chessboard, while $(k-1)$-dimensional faces of $\Delta_{m,n}$ are all
configurations of $k$ non-taking (non-attacking) rooks, meaning that two rooks are not allowed to be in the same row or the same column.

\medskip
The chessboard complex appears in different areas of mathematics and in many incarnations (as a coset complex of the symmetric group,
the matching complex in a complete bipartite graph, the complex of
all injective functions, etc.). Its topological properties (high connectivity and the structure of an orientable pseudomanifold)
have played a fundamental role in the proof of some deep results of topological combinatorics and discrete geometry (colored Tverberg theorems), see \cite{Z17} for a survey.

\begin{prop}(\cite{blvz, Z17, ZV92})\label{prop:chess}
The chessboard complex $\Delta_{m,n}$ is $(m-2)$-connected, $(m-1)$-dimensional simplicial complex for $n\geq 2m-1$.
\end{prop}

\medskip

\subsection*{Acknowledgements}

 R. \v Zivaljevi\' c was supported by the Science Fund of the Republic of Serbia, Grant No.\ 7744592, Integrability and Extremal Problems in Mechanics, Geometry and Combinatorics - MEGIC. Sections 2  and 4  are  supported by the Russian Science Foundation under grant  21-11-00040.

\end{document}